\input amstex
\documentstyle{amsppt}
\magnification=1200
\hcorrection{.25in}
\advance\vsize-.75in
\def\bF{\Bbb F}
\def\bZ{\Bbb Z}

\overfullrule=0pt
\topmatter
\title{A key equation and the computation of error values for 
codes from order domains}\endtitle
\leftheadtext{A key equation for codes from order domains}
\rightheadtext{A key equation for codes from order domains}
\author John B. Little \endauthor
\address Department of Mathematics and Computer Science,
         College of the Holy Cross,
         Worcester, MA 01610\endaddress
\email little\@mathcs.holycross.edu \endemail
\date  April 7, 2003 \enddate
\thanks{Research at MSRI is supported in part by NSF grant 
DMS-9810361.}\endthanks
\abstract
We study the computation of error values in the decoding of
codes constructed from order domains.  Our approach is based on a 
sort of analog of the {\it key equation} for decoding Reed-Solomon 
and BCH codes.  We identify a key equation for all codes
from order domains which have finitely-generated 
value semigroups; the field of fractions of the order domain
may have arbitrary transcendence degree, however.  We provide a 
natural interpretation of the construction using the theory of 
Macaulay's {\it inverse systems} and duality.  O'Sullivan's 
generalized Berlekamp-Massey-Sakata (BMS) decoding algorithm applies 
to the duals of suitable evaluation codes from these order domains.  
When the BMS algorithm does apply, we will show how it can be 
understood as a process for constructing a collection of solutions 
of our key equation. 
\endabstract
\endtopmatter

\document
\heading \S 1. Introduction \endheading

The theory of error control codes constructed using
ideas from algebraic geometry (including the geometric Goppa 
and related codes) has recently undergone a remarkable extension and 
simplification with the introduction of codes constructed from 
{\it order domains\/}.  Interestingly, this development has
been largely motivated by the structures utilized in the 
Berlekamp-Massey-Sakata decoding algorithm with Feng-Rao-Duursma 
majority voting for unknown syndromes.

We will review the definition of an order domain
in \S 2; for now we will simply say that the order domains form a 
class of rings having many of the same properties as the rings 
$R = \cup_{m=0}^\infty L(mQ)$ underlying the one-point geometric 
Goppa codes constructed from curves.  The general theory gives
a common framework for these codes, $n$-dimensional cyclic codes, as 
well as many other Goppa-type codes constructed from varieties of 
dimension $> 1$.  H\o holdt,  Pellikaan, and van Lint have given an 
exposition of order domains in \cite{HPL}, synthesizing work of 
many others in the coding theory community, and this is probably 
the best general reference for this topic.

More recently, Geil and Pellikaan (\cite{GP},\cite{Gei}) and 
O'Sullivan (\cite{OS1}) have studied the structure of order domains whose 
fields of fractions have arbitrary transcendence degree.  Moreover, 
O'Sullivan (\cite{OS2}) has shown that the Berlekamp-Massey-Sakata 
decoding algorithm (abbreviated as the BMS algorithm in the following)
and the Feng-Rao procedure extend in a natural way to a suitable 
class of codes in this much more general setting.

The decoding problem here can be divided into two parts:  
determination of the error locations, then determination of the 
corresponding error values.  When it applies, the BMS algorithm 
produces a Gr\"obner basis for what is known in the usual
terminology as the error-locator ideal corresponding
to the error vector, hence sufficient information to determine 
the error locations.  Here, we will consider the problem of 
determining the error values in conjunction with the BMS algorithm
or some other algorithm that determines the error locator ideal.  

For the Reed-Solomon codes (the simplest examples of codes from 
order domains, or geometric Goppa codes), the 
Berlekamp-Massey decoding algorithm (the precursor of BMS)
can be phrased as a method for solving a {\it key equation}.
For a Reed-Solomon code with minimum distance $d = 2t + 1$, 
the key equation has the form 
$$f S \equiv g \bmod \langle X^{2t}\rangle.\leqno(1.1)$$
Here ${ S}$ is a known univariate polynomial in $X$ constructed from 
the error syndromes, and $f,g$ are unknown polynomials in $X$.
If the error vector $e$ satisfies ${wt}(e) \le t$, there is
a unique solution $(f,g)$ with $\deg(f)\le t$, and $\deg(g)< \deg(f)$  
(up to a constant multiple).  The polynomial $f$ is known as 
the {\it error locator} because its roots give the {\it inverses} 
of the error locations; the polynomial $g$ is known as the 
{\it error evaluator} because the error values can be determined from 
values of $g$ at the roots of $f$, via the Forney formula.

O'Sullivan has introduced a generalization of this key 
equation for one-point geometric Goppa codes from curves 
in \cite{OS3} and shown that the BMS algorithm can be
modified to compute the analogs of the error-evaluator
polynomial together with error locators.  His
definitions make heavy use of the particular features
of the curve case, however.  For instance the 
objects corresponding to ${ S}$ and $g$ in (1.1) 
are differentials on the underlying curve.

Our main goals in this article are the following.
First, we wish to identify an analog of the key equation
(1.1) for codes from order domains.
We will only consider order domains whose value semigroups are
finitely generated.  In these cases, the ring $R$ can be 
presented as an affine algebra $R \cong \bF[X_1,\ldots,X_s]/I$,
where the ideal $I$ has a Gr\"obner basis of a very particular 
form (see \cite{GP} and \S 2 below).  
Although O'Sullivan has shown how more general order domains
arise naturally from valuations on function fields, it is not
clear to us how our approach applies to those examples. 
On the positive side, by basing all constructions 
on algebra in polynomial rings, all codes from these
order domains can be treated in a uniform way,
Second, we also propose to study the relation between the BMS 
algorithm and the process of solving this key equation in the 
cases where BMS is applicable.  Finally, we wish to show how 
solutions of our key equation can be used to determine error 
values and complete the decoding process.
 
Our key equation generalizes the key equation for 
$n$-dimensional cyclic codes studied by Chabanne and Norton in 
\cite{CN}.  Results on the algebraic background for their
construction appear in \cite{Norton1}.  See also \cite{Norton2} 
for connections with the more general problem of finding
shortest linear recurrences, and \cite{NS} for a generalization
giving a key equation for codes over commutative rings.  
In the present article, we will point out another natural 
interpretation of these ideas in the context of Macaulay's
{\it inverse systems} for ideals in a polynomial ring (see \cite{Mo}, 
\cite{EI}) and the theory of duality.

In spirit, our approach is also quite close to the treatment
of one-point geometric Goppa codes from curves 
by Heegard and Saints in \cite{HS}, in that we essentially 
treat all of our codes as (subcodes of) punctured $n$-dimensional 
cyclic codes. 

The present article is organized as follows.  In \S 2 we will
briefly review the definition of an order domain, evaluation
codes and dual evaluation codes.  We will also introduce some
standard examples.  
\S 3 contains a quick summary of the basics of Macaulay inverse
systems and duality for quotients of a polynomial ring by 
zero-dimensional ideals.  In \S 4 we introduce the key equation. 
We will also relate the BMS algorithm to the process of solving 
this equation.
\S 5 is devoted to a discussion of how the key equation
can be used to determine error values.  The 
major idea appears already for the case of $n$-dimensional
cyclic codes in \cite{CN}.  However, our results apply more 
generally and include a few improvements.  Finally, in 
\S 6 we present two detailed decoding examples using these methods. 

The author wishes to thank Mike O'Sullivan for comments on an earlier
draft of this paper.

\heading \S 2. Codes from Order Domains \endheading

In this section we will briefly recall the definition of
order domains and explain how they can be used to 
construct error control codes.  We will use the following 
formulation.

\bigskip
\noindent
{\bf (2.1) Definition.}  Let $R$ be a $\bF_q$-algebra and
let $(\Gamma,+,\succ)$ be a well-ordered semigroup.
An {\it order function} on $R$ is a surjective mapping 
$\rho : R \to \{-\infty\} \cup \Gamma$
satisfying:
\roster
\item $\rho(f) = -\infty \Leftrightarrow f = 0$,
\item $\rho(cf) = \rho(f)$ for all $f\in R$, all
$c \ne 0$ in $\bF_q$,
\item $\rho(f+g) \preceq \max_\succ \{\rho(f),\rho(g)\}$,
\item if $\rho(f) = \rho(g) \ne -\infty$, then there exists
$c \ne 0$ in $\bF_q$ such that $\rho(f) \prec \rho(f-cg)$,
\item $\rho(fg) = \rho(f) + \rho(g)$.
\endroster
\noindent
We call $\Gamma$ the {\it value semigroup} of $\rho$.

\bigskip
The terminology ``order function'' is supposed to 
suggest the existence of $\bF_q$-bases of $R$ 
whose elements have distinct $\rho$-values, and are
hence ordered by $\rho$.  This is a consequence of Axiom 4.
It is also possible to reindex the corresponding bases
by the natural numbers and define order functions in 
a different but equivalent way.  This is done, for
instance, in \cite{OS1} and \cite{OS2}.  

Axioms 1 and 5 in this definition imply that $R$ must be an
integral domain.  In the cases where the transcendence degree
of $R$ over $\bF_q$ is at least $2$, a ring $R$ with one
order function will have many others too.  For this reason
an {\it order domain} is formally defined as a pair $(R,\rho)$
where $R$ is an $\bF_q$-algebra and $\rho$ is an order function 
on $R$.  However, from now on,
we will only use one particular order function on $R$ 
at any one time.  Hence we will often omit it in refering 
to the order domain, and we will refer to $\Gamma$ as the 
value semigroup of $R$.  

From one point of view, order functions come from {\it valuations} 
on $K = QF(R)$.  As noted by O'Sullivan in \cite{OS1}, in fact 
$S = \{f/g : \rho(g) \succeq \rho(f)\}$ is a valuation ring of $K$.  
From now on, we will restrict our attention to the case that $\Gamma$
is a sub-semigroup of $\bZ_{\ge 0}^r$, for some $r \ge 1$, hence is 
{\it finitely generated}.  Without loss of generality, then, 
we may assume $r = tr.deg.{}_{\bF_q}(K)$.
To obtain a well-ordering on $\bZ_{\ge 0}^r$ we can fix a monomial 
order, $\succ$.

As noted in the introduction, order domains give a common 
generalization of several types of rings that have been 
used in the construction of codes.  For instance, the order domains used
in the construction of one-point geometric Goppa codes
are the following.  If $Y$ is a smooth projective curve
defined over $\bF_q$, and $Q$ is an $\bF_q$-rational point
on $Y$, then $R = \cup_{m=0}^\infty L(mQ)$ is an order domain.
$\Gamma$ is equal to the Weierstrass semigroup of $Y$ at $Q$ 
(the sub-semigroup of $\bZ_{\ge 0}$ consisting
of all pole orders of rational functions 
on $X$ with poles only at $Q$), and $\rho(f) = -v_Q(f)$, 
where $v_Q$ is the discrete valuation at $Q$ on the function field of $Y$.  
The polynomial ring $R = \bF_q[X_1,\ldots,X_r]$ is an order domain,
where $\Gamma = \bZ_{\ge 0}^r$, $\succ$ is a monomial
order, and $\rho(f)$ for $f \ne 0$ is defined by $\rho(f) = \alpha$
if $LT_{\succ}(f) = X^\alpha$.  These examples of order
domains feature in the construction of Reed-Muller and other
multidimensional cyclic codes.  Many other classes of examples
are considered in \cite{Gei} and \cite{GP}.

Geil and Pellikaan (see \cite{GP}) have proved a
characterization of order domains with finitely generated $\Gamma$,
which we will now review. 
In the following statement, $M$ is an $r\times s$ matrix with entries
in $\bZ_{\ge 0}$ with linearly independent rows.  
For $\alpha\in \bZ_{\ge 0}^s$ (written as 
a column vector), the matrix product $M\alpha$ is a vector in 
$\bZ_{\ge 0}^r$.  We will call this the $M$-weight of 
the monomial.  We write $\langle M\rangle$ for the subsemigroup 
of $\bZ_{\ge 0}^r$ generated by the {\it columns} of $M$,
ordered by any convenient monomial order $\succ$ on $\bZ_{\ge 0}^r$
(for instance the $lex$ order as in Robbiano's characterization
of monomial orders by weight matrices).
We will make use of the monomial orders $>_{M,\tau}$
on $\bF_q[X_1,\ldots,X_s]$ defined as follows:
$X^\alpha >_{M,\tau} X^\beta$ if 
$M\alpha \succ M\beta$, or if $M\alpha = M\beta$ and 
$X^\alpha >_\tau X^\beta$, where $\tau$ is another
monomial order used to break ties. 

\proclaim{(2.2) Theorem} (Geil-Pellikaan)
\roster
\item Let $\Gamma = \langle M\rangle \subset \bZ_{\ge 0}^r$ 
be a semigroup.  Let
$I \subset \bF_q[X_1,\ldots,X_s]$ be an ideal, and let $G$  
be the reduced Gr\"obner basis for $I$ with respect to a weight
order $>\ =\ >_{M,\tau}$ as above.  Suppose that 
every element of $G$ has exactly two monomials of highest 
$M$-weight in its support, and that the monomials 
in the complement of $LT_>(I)$ (the ``standard monomials''
or monomials in the ``footprint of the ideal'')
have distinct $M$-weights.  
Then $R = \bF_q[X_1,\ldots,X_s]/I$ 
is an order domain with value semigroup $\Gamma$  
and order function $\rho$ defined as follows:
Writing $f$ in $R$ as a linear combination
of the monomials in the complement of $LT_>(I)$,
$\rho(f) = \max_\succ\{M\beta: X^\beta \in supp(f)\}$.
\item  Every order domain with semigroup 
$\Gamma = \langle M\rangle$ has a presentation
$R \cong \bF_q[X_1,\ldots,X_s]/I$ such that the reduced
Gr\"obner basis of $I$ with respect to $>_{M,\tau}$ and 
the standard monomials have the form described in part (1).
\endroster
\endproclaim

In principle, this result gives a method to
construct the order domains with a given value semigroup 
$\Gamma$, as in following example.  

\bigskip
\noindent
{\bf (2.3) Example.}
Take $r = 2$, 
$\Gamma = \langle M \rangle \subset \bZ_{\ge 0}^2$,
ordered by $\succ$ the lexicographic order, where 
$$M = \pmatrix 0&1&3\cr
               2&1&0\cr\endpmatrix.$$
By the definition, the order function $\rho$ must be 
surjective, so there exist $x,y,z\in R$ with 
$\rho(x) =(0,2)$, $\rho(y) = (1,1)$, $\rho(z) = (3,0)$,
$R$ is generated by $x,y,z$, and $\rho(x^ay^bz^c)$ is
equal to the $M$-weight
$M (a,b,c)^t$ for all monomials $x^a y^b z^c$.
It follows that there is a surjective ring homomorphism
$$\phi : \bF_q[X,Y,Z] \to R,$$
where $\phi(X) = x$, $\phi(Y) = y$, and $\phi(Z) = z$.
We consider the monomial order $>_{M,lex}$ on $\bF_q[X,Y,Z]$.
It is easy to see that all $\bZ$-relations between $\rho(x)$,
$\rho(y)$, $\rho(z)$ are generated by $3\rho(x) + 2\rho(z) = 6\rho(y)$.
For Definition (2.1) to hold, we must have 
$\rho(x^2z^3 - cy^6) < \rho(x^2z^3)$ for some $c\ne 0$.
Hence $R \cong \bF_q[X,Y,Z]/I$, where $I = \langle F\rangle$
for some $F = X^2Z^3 - cY^6 + H(X,Y,Z)$, where every 
term in $H$ is less than $Y^6$ in the $>_{M,lex}$ order.  
The monomials in the complement of $\langle X^2Z^3\rangle$ 
have distinct $M$-weights and $G = \{F\}$ is 
a Gr\"obner basis for $I$ of the required form, so 
all such $R$ are order domains by Theorem (2.2).  Note that all are
deformations of the monomial algebra $\bF_q[u^2,uv,v^3]$.
Indeed, Theorem (2.2) can be reinterpreted as saying that 
the order domains with semigroup $\Gamma$ are all flat 
deformations of the monomial algebra $\bF_q[\Gamma]$.
This point of view is exploited in \cite{L} to construct 
order domains in the function fields of varieties
such as Grassmannians and flag varieties.

The most direct way to construct codes from an order domain
given by a particular presentation $R\cong\bF_q[X_1,\ldots,X_s]/I$ 
is to generalize Goppa's construction in the case of curves:

\bigskip
\noindent
{\bf (2.4) Construction of Codes.}
\roster
\item Let $X_R$ be the variety $V(I) \subset {\Bbb A}^s$ and 
let  
$$X_R(\bF_q) = \{P_1,\ldots,P_n\}$$
be the set of $\bF_q$-rational points on $X_R$.
\item Define an evaluation mapping 
$$\eqalign{ev : R &\to \bF_q^n\cr
                f &\mapsto (f(P_1),\ldots,f(P_n))\cr}$$
\item Let $V\subset R$ be any finite-dimensional vector
subspace.  Then the image $ev(V) \subseteq \bF_q^n$ 
will be a linear code in $\bF_q^n$.  One can also
consider the dual code $ev(V)^\perp$.
\item Of particular interest here are the codes
constructed as follows. Let $\Delta$ be the ordered basis of $R$
given by the monomials in the complement of $LT_>(I)$.
Note that this basis comes equiped with an ordering by $\rho$-value, or 
equivalently by the $M$-weights ordered by $\succ$ in $\bZ_{\ge 0}^r$.
Let $\ell \in {\Bbb N}$ and 
let $V_\ell$ be the span of the first
$\ell$ elements of the ordered basis $\Delta$. 
In this way, we obtain evaluation
codes $Ev_\ell = ev(V_\ell)$ and dual codes $C_\ell = Ev_\ell^\perp$
for all $\ell$.
\endroster

\bigskip
The BMS algorithm is specifically tailored for this last
class of codes.  If the $C_\ell$ codes are used 
to encode messages, then the $Ev_\ell$ codes describe
the parity checks and the syndromes used in the decoding
algorithm.

\heading \S 3. Preliminaries on Inverse Systems \endheading

A natural setting for our formulation of a key equation
for codes from order domains is the theory of inverse systems
of polynomial ideals originally introduced by Macaulay (\cite{Ma}).
There are several different versions of this theory. 
For modern versions using the language of differentiation
operators, see \cite{Mo} or \cite{EI}.  Here, we will summarize
a number of more or less well-known results, using an 
alternate formulation of the definitions that works in any 
characteristic.  A reference for this approach is \cite{North}.

Let $k$ be a field, let $S = k[X_1,\ldots,X_s]$ and
let $T$ be the formal power series ring $k[[X_1^{-1},\ldots,X_s^{-1}]]$
in the inverse variables.  $T$ is an $S$-module
under a mapping 
$$\eqalign{c : S \times T &\to T\cr
               (f,g) &\mapsto f\cdot g,\cr}$$
sometimes called contraction, defined as follows.  First, given monomials
$X^\alpha$ in $S$ and $X^{-\beta}$ in $T$, 
$X^\alpha \cdot X^{-\beta}$ is defined to be $X^{\alpha-\beta}$ 
if this is in $T$, and $0$ otherwise.  We then extend by linearity 
to define $c : S\times T \to T$.

Let $Hom_k(S,k)$ be the usual linear dual vector space.
It is a standard fact that the mapping 
$$\eqalign{\phi : Hom_k(S,k) &\to T\cr
                      \Lambda &\mapsto \sum_{\beta\in \bZ_{\ge 0}^s}
\Lambda(X^\beta) X^{-\beta}\cr}$$
is an isomorphism of $S$-modules, if we make $Hom_k(S,k)$
into an $S$-module in the usual way by defining
$(q\Lambda)(p) = \Lambda(qp)$ for all polynomials $p,q$ in $S$.
In explicit terms, the $k$-linear form on $S$ obtained from an 
element of $g \in T$ is mapping $\Lambda_g$ defined as follows.
For all $f\in S$, 
$$\Lambda_g(f) = (f\cdot g)_0,$$
where $(t)_0$ denotes the constant term in $t \in T$. 
In the following we will identify elements of $T$ with their
corresponding linear forms on $S$. 

For each ideal $I\subseteq R$, we can define the annihilator,
or {\it inverse system}, of $I$ in $T$ as
$$I^\perp = \{\Lambda \in T : \Lambda(p) = 0,\ \forall\  p \in I\}.$$
It is easy to check that $I^\perp$ is an $S$-submodule of $T$
under the module structure defined above.
Similarly, given an $S$-submodule $H \subseteq T$, we can define
$$H^\perp = \{p \in S : \Lambda(p) = 0,\ \forall\ \Lambda \in H\},$$
and $H^\perp$ is an ideal in $R$.

The key point in this theory is the following duality statement.

\proclaim{(3.1) Theorem}  The ideals of $R$ and the $S$-submodules
of $T$ are in inclusion-reversing
bijective correspondence via the constructions above,
and for all $I, H$ we have: 
$$(I^\perp)^\perp = I, \quad (H^\perp)^\perp = H.$$
\endproclaim

\noindent
See \cite{North} for a proof.

\bigskip
We will be interested in applying Theorem (3.1) when 
$I$ is the ideal of some finite set of points in 
the $n$-dimensional affine space over $k$ (e.g. when $k = \bF_q$ and $I$ is 
an error-locator ideal arising in decoding -- see \S 4 below).  

\proclaim{(3.2) Lemma}  Let
$$I = m_{P_1} \cap \cdots \cap m_{P_t},$$
where $m_{P_i}$ is the maximal ideal of 
$S$ corresponding to the point $P_i$, and $t\ge 1$.  The submodule
of $T$ corresponding to $I$ has the form
$$H = I^\perp = (m_{P_1})^\perp \oplus \cdots \oplus (m_{P_t})^\perp.$$
\endproclaim

\noindent
{\bf Proof.}
In Proposition 2.6 of \cite{Ger}, Geramita shows that
$(I\cap J)^\perp = I^\perp + J^\perp$ for any pair
of ideals.  The idea is that $I^\perp$ and $J^\perp$
can be constructed degree by degree, so the corresponding
statement from the linear algebra of finite-dimensional vector
spaces applies.  The equality $(I + J)^\perp = I^\perp \cap J^\perp$
also holds from linear algebra (and no finite-dimensionality
is needed).  The sum 
in the statement of the Lemma 
is a direct sum since $m_{P_i} + \cap_{j\ne i} m_{P_j} = S$,
hence $(m_{P_i})^\perp \cap \Sigma_{j\ne i} (m_{P_j})^\perp = \{0\}$.
$\square$

\bigskip
We can also give a concrete description of the elements
of $(m_P)^\perp$.

\proclaim{(3.3) Proposition} Let $P = (a_1,\ldots,a_s) \in {\Bbb A}^s$
over $k$, and let $L_i$ be the coordinate hyperplane $X_i = a_i$ containing 
$P$.
\roster
\item $(m_P)^\perp$ is the cyclic $S$-submodule of $T$ generated by
$$h_P = \sum_{u\in \bZ_{\ge 0}^s} P^u X^{-u},$$
where if $u = (u_1,\ldots,u_s)$, $P^u$ denotes the product
$a_1^{u_1} \cdots a_s^{u_s}$ ($X^u$ evaluated at $P$).
\item $f\cdot h_P = f(P)h_P$ for all $f\in S$,
and the submodule $(m_P)^\perp$ is a one-dimensional vector space over $k$.
\item Let $I_{L_i}$ be the ideal $\langle X_i - a_i\rangle$ in 
$S$ (the ideal of $L_i$).  
Then $(I_{L_i})^\perp$ is the submodule of $T$ generated by
$$h_{L_i} = \sum_{j=0}^\infty a_i^j X_i^{-j}.$$
\item In $T$, we have 
$$h_P = \prod_{i=1}^s h_{L_i}.$$
\endroster
\endproclaim

\bigskip
\noindent
{\bf Proof.}  (1) First, if $f \in m_P$, and $g\in S$ is arbitrary then 
$$\Lambda_{g\cdot h_P}(f) = (f\cdot (g\cdot h_P))_0 = ((fg)\cdot h_P)_0 = f(P)g(P) = 0.$$ 
Hence the $S$-submodule $\langle h_P\rangle$ is contained in 
$(m_P)^\perp$.  Conversely, if $h \in (m_P)^\perp$, then for all $f\in m_P$,
$$0 = \Lambda_h(f) = (f\cdot h)_0.$$ 
An easy calculation using all $f$ of the form $f= x^\beta - a^\beta \in m_P$ shows that  
$h = c h_P$ for some constant $c$.  Hence $(m_P)^\perp = \langle h_P\rangle$. 

(2) The second claim follows by a direct computation of the contraction
product $f\cdot h_p$.

(3) Let $f\in I_{L_i}$ (so $f$ vanishes at all points of the 
hyperplane $L_i$), and let $g\in S$ be arbitrary.  Then
$$\eqalign{\Lambda_{g\cdot h_{L_i}} (f) 
&= (f\cdot (g\cdot h_{L_i}))_0\cr
&= ((fg)\cdot h_{L_i})_0\cr 
&= f(0,\ldots,0,a_i,0,\ldots,0) g(0,\ldots,0,a_i,0,\ldots,0)\cr
&= 0,\cr}$$
since the only nonzero terms in the product $((fg)\cdot h_{L_i})$
come from monomials in $fg$ containing only the variable $X_i$.
Hence $\langle h_{L_i}\rangle \subset T$ is contained
in $I_{L_i}^\perp$.  Then we show the other inclusion as in the 
proof of (1).

(4) We have $m_P = I_{L_1} + \cdots + I_{L_s}$.  Hence
$(m_P)^\perp = (I_{L_1})^\perp \cap \cdots \cap (I_{L_s})^\perp$, and 
the claim follows.  We note that a more explicit form of this
equation can be derived by the formal geometric series summation
formula:
$$h_P = \sum_{u\in \bZ_{\ge 0}^s} P^uX^{-u} = \prod_{i=1}^s {1\over 1-a_i/X_i} = 
\prod_{i=1}^s h_{L_i}.\ \square$$ 

\bigskip

Finally, we note that both the polynomial ring $S$ and the 
formal power series ring $T$ can be viewed as subrings of the 
field of formal Laurent series in the inverse variables, 
$$K = k((X_1^{-1}, \ldots, X_s^{-1})),$$
which is the field of fractions of $T$.
Hence there is a natural interpretation of the (full) product
$fg$ for $f \in S$ and $g \in T$ as an element of $K$.  
The contraction product
$f\cdot g$ can be understood as a projection of $fg$ into
$T \subset K$ (image under the linear
projection with kernel spanned by all monomials
not in $T$).  In the sequel, we will also need to make use of the 
projection of $fg$ into $S_+ = \langle X_1,\ldots, X_s\rangle \subset S\subset K$
under the linear projection with kernel spanned by all monomials not in $S_+$.
We will denote this by $(fg)_+$.  Hence $(fg)_+$ gives the 
sum of all terms in $fg$ with all exponents nonnegative and some
exponent strictly positive, while $f\cdot g$ gives the 
sum of all terms in $fg$ with nonpositive exponents.  Any
``mixed terms'' in $fg$ (i.e. those terms with some positive and 
some negative exponents) will be irrelevant in our 
applications.  We will use the following fact.

\proclaim{(3.4) Proposition}  Let $f\in k[X_i]$
be a univariate polynomial satisfying $f(P) = 0$.  Then
$$(fh_P)_+ = X_i g,$$
where $g(P) = f'(P)$ (formal derivative).
\endproclaim

\noindent
{\bf Proof.}  This follows by a direct computation using (3.3). $\square$

\heading \S 4. The Key Equation and its Relation to the
BMS Algorithm\endheading

In this section, we will introduce our key equation
for codes from order domains and relate it to the 
Berlekamp-Massey-Sakata decoding algorithm.
Let $C$ be one of the codes $C = ev(V)$ or $ev(V)^\perp$ 
constructed from an order domain 
$R \cong \bF_q[X_1,\ldots,X_s]/I$ as in \S 2 above.  
Consider an error vector $e \in \bF_q^n$ 
(where entries are indexed by the elements of the set $X_R(\bF_q)$).  
In the usual terminology, the {\it error-locator ideal} corresponding
to $e$ is the ideal $I_e \subset \bF_q[X_1,\ldots,X_s]$
defining the set of error locations: 
$$I_e = \{ f \in \bF_q[X_1,\ldots,X_s] : f(P) = 0,\ \forall\ P\ s.t.\ e_P \ne 0\}.$$
(Since $I_e \supset I$, one could also consider the ideal corresponding  
to $I_e$ in $R$.   However, 
following the general philosophy of Heegard and Saints in \cite{HS}, 
we will find it more convenient to work with $I_e$ as an ideal
in the polynomial ring.)

We will also use a slightly different notation and terminology
in the following because we want to make a systematic use
of the observation that this ideal depends only on the 
support of $e$, not on the error values.  Indeed, many different
error vectors yield the same ideal defining
the error locations.  For this reason we will introduce
${\Cal E} = \{P : e_P \ne 0\}$, and refer to the error-locator
ideal for any $e$ with $supp(e) = {\Cal E}$ as $I_{\Cal E}$.

For each monomial $X^u \in \bF[X_1,\ldots,X_s]$, we let
$$E_u = \langle e, ev(X^u)\rangle = \sum_{P\in X_R(\bF_q)} e_P P^u
\leqno(4.1)$$
be the corresponding syndrome of the error vector.  (As in (3.3),
$P^u$ is shorthand notation for the evaluation of the monomial
$X^u$ at $P$.)
  
In the practical decoding situation, of course, 
for a code $C = ev(V)^\perp$ where $V$ is a subspace
of $R$ spanned by some set of monomials, 
only the $E_u$ for the $X^u$ in a basis of $V$
are initially known from the received word.  

In addition, the elements of the ideal
$I + \langle X_1^q-X_1,\ldots,X_s^q - X_s\rangle$
defining the set $X_R(\bF_q)$ give relations between the 
$E_u$.  Indeed, the $E_u$ for $u$ in the ordered basis $\Delta$ for $R$ 
with all components $\le q - 1$ determine all the others, 
and these syndromes still satisfy additional relations.
Thus the $E_u$ are, in a sense, highly redundant.  

To package the syndromes into a single algebraic
object, we define the {\it syndrome series}
$${\Cal S}_e = \sum_{u\in \bZ_{\ge 0}^s} E_u X^{-u}$$
in the formal power series ring
$T = \bF_q[[X_1^{-1},\ldots,X_s^{-1}]]$.  (This depends
both on the set of error locations ${\Cal E}$ and on 
the error values.)
Chabanne and Norton considered the same type of 
expression in \cite{CN} for $n$-dimensional cyclic codes.
As in \S 3, we have a natural interpretation for 
${\Cal S}_e$ as an element of the dual space of the ring
$S = \bF_q[X_1,\ldots,X_s]$.  

A fundamental tool in our considerations will
be the following expression 
for the syndrome series ${\Cal S}_e$.  We substitute 
from (4.1) for the syndrome $E_u$ and change the order of 
summation to obtain:
$$\eqalign{
{\Cal S}_e &= \sum_{u\in \bZ_{\ge 0}^n} E_u X^{-u}\cr
&= \sum_{u\in \bZ_{\ge 0}^n} \sum_{P\in X_R(\bF_q)} e_P P^u X^{-u}\cr
&= \sum_{P\in X_R(\bF_q)} e_P \sum_{u\in \bZ_{\ge 0}^n} P^u X^{-u}\cr
&= \sum_{P\in X_R(\bF_q)} e_P h_P,\cr}\leqno(4.2)$$
where $h_P$ is the generator of $(m_P)^\perp$ from (3.3).
The sum in (4.2), taking the terms with $e_P\ne 0$, 
gives the decomposition of ${\Cal S}_e$ in the direct
sum expression for $I_{\Cal E}^\perp$ as in (3.2).

The following result
is well-known in a sense; it is a translation of the 
standard fact that error-locators give linear recurrences
on the syndromes.  But to our knowledge, this connection has not
been considered from exactly our point of view in
this generality (see \cite{AD} for a special case).

\bigskip
\proclaim{(4.3) Theorem} With all notation as above,  
\roster
\item  $f\in I_{\Cal E}$ if and only if $f\cdot {\Cal S}_e=0$ for
all error vectors $e$ with $supp(e) = {\Cal E}$.
\item For each $e$ with $supp(e) = {\Cal E}$, 
$I_{\Cal E} = \langle {\Cal S}_e\rangle^\perp$ in the duality from
Theorem (3.1).
\item If $e,e'$ are two error vectors with the same support, then
$\langle {\Cal S}_e \rangle = \langle {\Cal S}_{e'}\rangle$ 
as submodules of $T$.
\endroster
\endproclaim

\bigskip
\noindent
{\bf Proof.}  For (1), we start from the expression for 
${\Cal S}_e$ from (4.2).  Then by (3.3), we have
$$f\cdot {\Cal S}_e = \sum_{P\in {\Cal E}} e_P (f\cdot h_P) = 
\sum_{P\in {\Cal E}} e_P f(P) h_P.$$
If $f \in I_{\Cal E}$, then clearly 
$f \cdot {\Cal S}_e = 0$ for all choices of error values
$e_P$.  Conversely,  if $f\cdot {\Cal S}_e = 0$ for all $e$ with
$supp(e) = {\Cal E}$, then $f(P) = 0$ for all $P \in {\Cal E}$,
so $f\in I_{\Cal E}$.

Claim (2) follows from (1).  

The perhaps surprising claim (3) is a consequence of (2).  Another way
to prove (3) is to note that there exist $g\in R$ 
such that $g(P)e_P = e'_P$ for all $P \in {\Cal E}$.  We have
$$g\cdot {\Cal S}_e = \sum_{P\in {\Cal E}} e_P (g\cdot h_P)
=\sum_{P\in {\Cal E}} e_P g(P) h_P = 
\sum_{P\in {\Cal E}} e'_P h_P = {\Cal S}_{e'}.$$
Hence $\langle {\Cal S}_{e'}\rangle \subseteq \langle {\Cal S}_e\rangle$.
Reversing the roles of $e$ and $e'$, we get the other inclusion
as well, and (3) follows.
$\square$

\bigskip
The following explicit expression for 
the terms in $f\cdot {\Cal S}_e$ is 
also useful. Let $f = \sum_m f_m X^m \in S$.  Then 
$$\eqalign{
f\cdot {\Cal S}_e &= (\sum_m f_m X^m)\cdot (\sum_{u\in \bZ_{\ge 0}^s} E_u X^{-u})\cr
          &= \sum_{r\in \bZ_{\ge 0}^s} (\sum_m f_m E_{m+r}) X^{-r}. \cr}\leqno(4.4)$$
Hence $f\cdot {\Cal S}_e = 0 \Leftrightarrow \sum_m f_m E_{m+r} = 0$
for all $r \ge 0$.

The equation $f\cdot {\Cal S} = 0$ from (1) in (4.3) is the prototype, 
so to speak, for our generalizations 
of the key equation to all codes from order domains, and we will
refer to it as the key equation in the following.  
It also naturally generalizes all the various key equations that have 
been developed in special cases, as we will demonstrate
shortly.  Before proceeding with that, however, we wish to 
make several comments about the form of this equation.

Comparing the equation $f\cdot {\Cal S}_e = 0$
with the familiar form (1.1), several differences 
may be apparent.  First, note that the syndrome 
series ${\Cal S}_e$ will not be entirely
known from the received word in the decoding situation.  
The same is true in the Reed-Solomon case, of course.
The polynomial $S$ in the congruence in (1.1) involves only 
the known syndromes, and (1.1) is derived by accounting for the 
other terms in the full syndrome series.  With a
truncation of ${\Cal S}_e$ in our situation we would obtain
a similar type of congruence (see the discussion following
(4.14) below, for instance).

It is apparently rare, however, that
the portion of ${\Cal S}_e$ known from the received word
suffices for decoding up to half the minimum distance
of the code.  As first noted for the one-point geometric Goppa codes
from curves, it is often the case that additional syndromes
(or other extra information about the error) must be determined 
in order to exploit the code's full error correcting capacity.  
For this reason, even though we have not made
any hypotheses so far on how our code was constructed (i.e. on
how the vector subspace $V \subset R$ was chosen), the key equation
will be most useful in the case that $C$ is one of the
codes $C_\ell = Ev_\ell^\perp$ defined in \S 2, for which
the Feng-Rao majority voting process for unknown syndromes and
the generalized BMS algorithm are applicable.

Another difference is that there is no apparent analog of the 
error-evaluator polynomial $g$ from 
(1.1) in the equation in $f\cdot {\Cal S}_e = 0$.  In \S 5, we will see 
that the way to obtain error evaluators in this situation is to consider 
the ``purely positive parts'' $(f{\Cal S}_e)_+$ for certain
solutions of our key equation.  

We now turn to several examples that show how our key 
equation relates to several special cases that have
appeared in the literature.  

\bigskip
\noindent
{\bf (4.5) Example.}  We begin by providing more detail 
on the precise relation between (4.3), part (1) in the case of a
Reed-Solomon code and the usual key equation from (1.1).  
These codes are constructed from the order domain $R = \bF_q[X]$ 
(where $\Gamma = \bZ_{\ge 0}$ and $\rho$ is the degree mapping),
according to (4.4).  The key equation (1.1) applies to
the code $Ev_\ell = ev(V_\ell)$,
where $V_\ell = Span\{1,X,X^2,\ldots,X^{\ell-1}\}$, 
and the evaluation takes place at all $\bF_q$-rational points
on the affine line, omitting $0$.

For the $Ev_\ell$ Reed-Solomon codes, the known syndromes
are $E_1,\ldots,E_{d-1}$, and $S$ is the syndrome polynomial:
$$S = E_1 + E_2 X + \cdots + E_{d-1} X^{d-2}.$$
In the special solution 
$(f,g)$ of (1.1) used for decoding, 
$$f = \prod_{i=1}^{|{\Cal E}|} (1 - \alpha^{e_i}x),$$
where $\alpha^{-e_i}$ are the error locations.  Moreover,
$$g = \sum_{i=1}^{|{\Cal E}|} e_i \alpha^i \prod_{j\ne i} (1 - \alpha^{e_j}x).$$
If (1.1) is written as an equation
$$fS = g + x^{2t} h,$$
then $h$ is another polynomial of degree $|{\Cal E}| - 1$ sometimes
called the {\it error coevaluator}:
$$h = \sum_{i=1}^{|{\Cal E}|} e_i \alpha^{(2t+1)i} \prod_{j\ne i} (1 - \alpha^{e_j}x).$$
Either $g$ or $h$ can be used to solve for the error 
values $e_i$ once the roots of $f$ are determined.

Our key equation in this case is closely related, but
not precisely the same.  The natural way to apply (4.3)
here is to the dual code $C_\ell = Ev_\ell^\perp$.
Our prototype key equation $f\cdot {\Cal S}_e = 0$ 
uses the full syndrome series, but of course, we
could also consider the truncation of ${\Cal S}_e$
using only the known syndromes $E_0,\ldots,E_{\ell-1}$
and obtain a congruence close in form to (1.1).

Starting from (4.4)
and using the formal geometric
series summation formula as in (3.3) part (4), we can write:
$$\eqalign{
{\Cal S}_e &= \sum_{P\in {\Cal E}} e_P h_P\cr
           &= \sum_{P\in {\Cal E}} e_P \sum_{u\ge 0} P^u X^{-u}\cr
           &= \sum_{P\in {\Cal E}} e_P {1\over 1 - P/X}\cr
           &= X {\sum_{P\in {\Cal E}} e_P 
\prod_{Q\in {\Cal E}, Q\ne P}(X - Q) \over \prod_{P \in {\Cal E}} (X - P)}\cr}$$
Hence, in this formulation, ${\Cal S}_e = Xq/p$, where
$p$ is the generator of the (actual) error locator ideal.
By considering the truncated form of $f\cdot {\Cal S}_e = 0$,
it can be seen that our $q$ is actually the 
analog of the error coevaluator as above. Moreover if $f = p$,
then $(p {\Cal S}_e)_+ = Xq$ gives the error (co)evaluator.  There
are no ``mixed terms'' in the products $f{\Cal S}_e$ in this one-variable
situation.

\bigskip
\noindent
{\bf (4.6) Example.} 
The key equation for $s$-dimensional cyclic codes introduced
by Chabanne and Norton in \cite{CN} has the form:
$$\sigma {\Cal S}_e = \left(\prod_{i=1}^s X_i\right) g,\leqno(4.7)$$
where 
$$\sigma = \prod_{i=1}^s \sigma_i(X_i),$$
and $\sigma_i$ is the univariate generator of the elimination 
ideal $I_{\Cal E} \cap \bF_q[X_i]$.  Our version of the Reed-Solomon key 
equation from (4.5) is a special case of (4.7).  Moreover, 
(4.7) is clearly the special case of (4.3), part (1) for these codes
where $f = \sigma$ is the particular error locator polynomial
$\prod_{i=1}^s \sigma_i(X_i) \in I_{\Cal E}$.  For this 
special choice of error locator, $\sigma \cdot {\Cal S}_e = 0$, and
$(\sigma{\Cal S}_e)_+ = \left(\prod_{i=1}^s X_i\right) g$ for some 
polynomial $g$.  This last claim can be established
using (4.4).  We see that ${\Cal S}_e$
can be written as
$${\Cal S}_e = \sum_P e_P h_P = \left(\prod_{i=1}^s X_i\right) 
\sum_P e_P {1\over \prod_{i=1}^s (X_i - X_i(P))}$$
and the product $\sigma {\Cal S}_e = (\sigma {\Cal S}_e)_+$ 
reduces to a polynomial (again, there are no ``mixed terms'').  

In order to use (4.7) for decoding, Chabanne and Norton propose
iterated applications of the one-variable Berlekamp-Massey algorithm
to find the factors of the product $\sigma$ one at a time.  
In \S 5 and \S 6 we will see that the more general BMS algorithm
gives additional flexibility for decoding these codes,
although the equation (4.7) will still lead most directly to 
determination of the error values.

\bigskip
\noindent
{\bf (4.8) Example.}  We now turn to the key equation for
one-point geometric Goppa codes introduced by O'Sullivan in 
\cite{OS3}.   Let ${\Cal X}$ be a smooth curve over $\bF_q$ of genus
$g$, and consider one-point codes constructed from 
$R = \cup_{m=0}^\infty L(mQ)$ for some point $Q \in {\Cal X}(\bF_q)$, 
O'Sullivan's key equation has the form:
$$f \omega_e = \phi.\leqno(4.9)$$
Here $\omega_e$ is the syndrome differential, which 
can be expressed as
$$\omega_e = \sum_{P\in {\Cal X}(\bF_q)} e_P \omega_{P,Q},$$
where $\omega_{P,Q}$ is the differential of the third kind on $Y$
with simple poles at $P$ and $Q$, no other poles, and residues
$$res_P(\omega_{P,Q}) = 1,\quad res_Q(\omega_{P,Q}) = -1.$$
For any $f \in R$, we have
$$res_Q(f\omega_e) = \sum_P e_P f(P),$$
the syndrome of $e$ corresponding to $f$.  (We only defined syndromes
for monomials above; taking a presentation $R = \bF_q[X_1,\ldots,X_s]/I$, however, any $f\in R$ can be expressed as a linear combination of monomials
and the syndrome of $f$ is defined accordingly.)
The right-hand side of (4.9) is also a differential.  In this 
situation, (4.9) furnishes a key equation in the following 
sense: $f$ is an error locator (i.e. $f$ is in the ideal of $R$
corresponding to $I_{\Cal E}$) if 
and only if $\phi$ has poles only at $Q$.  
In the special case that $(2g - 2)Q$ is a canonical divisor
(the divisor of zeroes of some differential of the first kind
$\omega_0$ on ${\Cal X}$), (4.9) can be replaced by the equivalent equation
$$f o_e = g,\leqno(4.10)$$
where $o_e = \omega_e/\omega_0$ and $g = \phi/\omega_0$ are rational
functions on ${\Cal X}$.  Since $\omega_0$ is zero only at $Q$,
the key equation is now that $f$ is an error locator if and only
if (4.9) is satisfied for some $g \in R$.

For instance, when ${\Cal X}$ is a smooth plane curve
$V(F)$ over $\bF_q$ defined by $F \in \bF_q[X,Y]$, 
with a single point $Q$ at infinity, then it is true that 
$(2g-2)Q$ is canonical.  O'Sullivan shows in Example 4.2 of 
\cite{OS3} (using a slightly different notation) that 
$$o_e = \sum_{P\in {\Cal X}(\bF_q)} e_P H_P,\leqno(4.11)$$
where if $P = (a,b)$, then $H_P = {F(a,Y)\over (X-a)(Y-b)}$.
This is a function with a pole of order $1$ at $P$, a pole
of order $2g - 1$ at $Q$, and no other poles.  

To relate this to our approach, note that we may assume from
the start that $Q=(0:1:0)$ and that $F$ is taken in the 
form from Theorem (2.2), that is 
$$F(X,Y) = X^\beta - cY^\alpha + G(X,Y)$$
for some relatively prime $\alpha < \beta$ generating the 
value semigroup at $Q$.  Every term in $G$ has 
$(\alpha,\beta)$-weight less than $\alpha\beta$.  

Then we can proceed as in Example (4.3) of \cite{OS3} to 
relate $H_P$ to an element of $T = \bF_q[[X^{-1},Y^{-1}]]$.
First we rearrange to obtain
$$\eqalign{H_P &= {F(a,Y)\over (X-a)(Y-b)}\cr
&={a^\beta - cY^\alpha + G(a,Y)\over (X-a)(Y-b)}\cr
&={(a^\beta - X^\beta) + F(X,Y) + (G(a,Y) - G(X,Y))\over (X-a)(Y-b)}\cr}$$
The $F(X,Y)$ term in the numerator does not depend on $P$.  We 
can collect those terms in the sum (4.11) and factor out the 
$F(X,Y)$.   We will see shortly that those
terms can in fact be ignored.  The $G(a,Y) - G(X,Y)$ in the numerator
furnish terms that go into the error evaluator $g$ here.
The remaining portion is
$${-(X^\beta - a^\beta)\over (X-a)(X-b)} = 
-{X^{\beta-1}\over Y}\sum_{i=0}^\beta \sum_{j=0}^\infty {a^ib^j\over X^iY^j}.$$
The sum here looks very much like that defining
our $h_P$ from (3.3), except that 
it only extends over the monomials in complement of $\langle LT(F)\rangle$.
Call this last sum $h_P'$.  As noted before the full series $h_P$ (and
consequently ${\Cal S}$) are redundant. 
For example, every ideal contained in $m_P$
(for instance the ideal $I = \langle F\rangle$ defining the curve), produces 
relations between the coefficients.  From the duality theorem (3.1),
we have that $I \subset m_P$ implies $(m_P)^\perp \subset I^\perp$,
so $F\cdot h_P = 0$.  

The relation $F\cdot h_P = 0$
says in particular that the terms in $h_P'$ are
sufficient to determine the whole series $h_P$.  Indeed, 
we have 
$$\eqalign{h_P &= h_P' + {(cY^\alpha - G)\over X^\beta}\cdot h_P' + 
\left({(cY^\alpha - G)\over X^\beta}\right)^2\cdot h_P' + \cdots\cr
&= \left({1\over 1 - {(cY^\alpha - G)\over X^\beta}}\right)\cdot h_P'\cr
&= \left({X^\beta \over F}\right)\cdot h_P'\cr}$$
It follows that O'Sullivan's key equation and ours are equivalent.
\bigskip

We now turn to the precise relation between solutions of 
our key equation and the polynomials generated by
steps of the BMS decoding algorithm applied to the $C_\ell = Ev_\ell^\perp$
codes from order domains $R$.  We will see that the steps
of the BMS algorithm systematically produce successively better 
approximations to solutions of $f\cdot {\Cal S}_e = 0$, so that in effect, 
{\it the BMS algorithm is a method for solving the key equation 
for these codes}.  
In addition to \cite{OS3} cited previously, a similar interpretation 
of the Berlekamp-Massey algorithm in the Reed-Solomon case (and related 
cases) was developed by Fitzpatrick in \cite{F} (see also
\cite{CLO}, Chapter 9, \S 4).  

We recall 
the key features of O'Sullivan's presentation of BMS. 
For our purposes, it will suffice to consider
the {\it ``Basic Algorithm''} from \S 3 of \cite{OS2}, in which
all needed syndromes are assumed known and no sharp
stopping criteria are identified. The 
{\it syndrome mapping} corresponding to the error vector $e$ is 
$$\eqalign{Syn_e : R &\to \bF_q\cr
                   f &\mapsto \sum_{P\in {\Cal E}} e_P f(P),\cr}$$
where as above ${\Cal E}$ is the set of error locations.
The same reasoning used in the proof of our Theorem (4.2) shows
$$f \in I_{\Cal E} \Leftrightarrow Syn_e(fg) = 0, \forall g \in R.\leqno(4.12)$$

From Definition (2.1) and Geil and Pellikaan's presentation
theorem (2.2), we have an ordered monomial basis of $R$:
$$\Delta = \{X^{\alpha(j)} : j \in {\Bbb N}\},$$  
whose elements have distinct $\rho$-values.  As in the 
construction of the $Ev_\ell$ codes, we write 
$V_\ell = Span\{1=X^{\alpha(1)},\ldots,X^{\alpha(\ell)}\}$.
The $V_\ell$ exhaust $R$, so for $f\ne 0\in R$, we may define
$$o(f) = \min \{\ell : f \in V_\ell\},$$
and (for instance) $o(0) = -1$.
Indeed, all properties of order domains can be
restated in terms of $o$, and O'Sullivan uses this function 
rather than $\rho$ in \cite{OS1} and \cite{OS2}.
In particular the semigroup $\Gamma$ in our presentation 
carries over to a (nonstandard) semigroup structure
on ${\Bbb N}$ defined by the addition operation
$$i\oplus j = k \Leftrightarrow o(X^{\alpha(i)}X^{\alpha(j)})=k.$$ 

Given $f\in R$, one defines
$$\eqalign{
span(f) &= \min\{\ell : \exists g \in V_\ell\ s.t.\ Syn_e(fg) \ne 0\}\cr
fail(f) &= o(f)\oplus span(f).\cr}$$
When $f\in I_{\Cal E}$, $span(f) = fail(f) = \infty$.

The BMS algorithm, then, is an iterative process which produces
a Gr\"obner basis for $I_{\Cal E}$ with respect to the monomial
order $>\ =\ >_{M,\tau}$ in (2.2).  The strategy is to maintain 
data structures for all $m\ge 1$ as follows.  The $\Delta_m$
are an increasing sequence of sets of monomials, converging
to the monomial basis for $I_{\Cal E}$ as $m \to \infty$. 
$\delta_m$ is the set of maximal elements of $\Delta_m$ 
with respect to $>$ (the ``interior corners of the footprint'').
Similarly, we consider $\Sigma_m = \bZ_{\ge 0}^s\ \backslash\ \Delta_m$,
and $\sigma_m$, the set of minimal elements of $\Sigma_m$ (the ``exterior
corners'').  
For sufficiently large $m$, the elements of $\sigma_m$ will
be the leading terms of the elements of the Gr\"obner basis of $I_{\Cal E}$,
and $\Sigma_m$ will the be set of monomials in $LT_>(I_{\Cal E})$.

For each $m$, the algorithm also produces collections of polynomials
$F_m = \{f_m(s) : s \in \sigma_m\}$ and $G_m = \{g_m(c) : c \in \delta_m\}$
satisfying:
$$o(f_m(s)) = s, \quad fail(f_m(s)) > m$$
and
$$span(g_m(c)) = c, \quad fail(g_m(c)) \le m.$$  
In the limit as $m\to \infty$, by (4.12), the $F_m$ yield
the Gr\"obner basis for $I_{\Cal E}$.

We record the following simple observation.

\proclaim{(4.13) Proposition}  With all notation as above,
suppose $f\in R$ satisfies $o(f) = s$, $fail(f) > m$.  Then
$$f\cdot {\Cal S}_e \equiv 0 \bmod W_{s,m},$$
where $W_{s,m}$ is the $\bF_q$-vector subspace of the formal
power series ring $T$ spanned by the 
$X^{-\alpha(j)}$ such that $s \oplus j > m$.
\endproclaim

\noindent
{\bf Proof.}  By the definition, $fail(f) > m$ means
that $Syn_e(f X^{\alpha(k)}) = 0$ for all $k$ with 
$o(f)\oplus k \le m$.  By the definitions 
of ${\Cal S}_e$ and the contraction product,
$Syn_e(fX^{\alpha(k)})$ is exactly the coefficient of
$X^{-\alpha(k)}$ in $f\cdot {\Cal S}_e$.  
$\square$

\bigskip
The subspace $W_{s,m}$ in (4.13) depends on $s=o(f)$.  In our
situation, though, note that if $s' = \max\{o(f) : f \in F_m\}$, 
then (4.13) implies
$$f \cdot {\Cal S}_e \equiv 0 \bmod W_{s',m}\leqno(4.14)$$
for all $f = f_m(s)$ in $F_m$.  Moreover, only finitely
many terms from ${\Cal S}_e$ enter into any one of these
congruences, so (4.14) is, in effect,
a sort of general analog of (1.1).

The $f_m(s)$ from $F_m$ can be understood as approximate
solutions of key equation (where the goodness of the approximation
is determined by the subspaces $W_{s',m}$,
a decreasing chain, tending to $\{0\}$ in $T$, as $m\to\infty$).  
The BMS algorithm thus systematically constructs better and better 
approximations to solutions of the key
equation.  O'Sullivan's stopping criteria (\cite{OS2})
show when further steps of the algorithm make no
changes.  Also note that the Feng-Rao theorem shows
that any additional syndromes needed for this can
be determined by the majority-voting process when
$wt(e) \le \lfloor{d_{FR}(C_\ell) - 1\over 2}\rfloor$.

We conclude this section by noting that O'Sullivan
has also shown in \cite{OS3} that, for codes from curves, 
the BMS algorithm can be slightly modified to compute error 
locators and error evaluators simultaneously in the situation 
studied in Example (4.7).  The same is almost certainly true
in our general setting, although we have not worked out all 
the details.  One reason we have not done so
is that it is not clear that all of the purely 
positive parts $(f{\Cal S}_e)_+$
for $f\in I_{\Cal E}$ are directly useful for determining
error values.  That seems to be true only for 
special $f\in I_{\Cal E}$ (in particular, for the univariate 
polynomials in the elimination ideals $I_{\Cal E}\cap \bF_q[X_i]$).  
In the practical decoding situation, once the BMS algorithm
is executed, the next step would be to solve a 
system of polynomial equations to determine
the error locations, i.e. to find the variety $V(I_{\Cal E})$
using the computed Gr\"obner basis for $I_{\Cal E}$.  Many of 
the same techniques useful for that process can efficiently
produce the needed univariate polynomials as a byproduct.
Hence we will not consider the sort of modification of
BMS proposed in \cite{OS3}.

\heading \S 5. Determination of Error Values \endheading

\noindent

In this section, we will see how solutions $f$ of the 
key equation (4.3), part (1) can be used to determine error values.
The method is the same as that presented in \cite{CN};
our proofs are significantly simplified by the use of the 
formalism from \S 3.

We will begin with some general results concerning
the polynomials $(f{\Cal S})_+$ for univariate $f\in I_{\Cal E}$. 
First we consider a simple special case.  Let 
${\Cal E} =\{P_1,\ldots,P_t\} = supp(e)$ for the error vector
$e$.  We will say that ${\Cal E}$ is in {\it general
position with respect to $X_i$} if the $X_i$-coordinates
of the $P_j$ are distinct.  

\proclaim{(5.1) Proposition} Let $e$ be an error
vector such that ${\Cal E}$ is
in general position with respect to $X_i$.   
Let $f$ be the monic generator of the elimination
ideal $I_{\Cal E} \cap \bF_q[X_i]$, then $(f{\Cal S}_e)_+ = X_i g$
for some $g \in \bF_q[X_i]$.  Moreover, if $P$
is any one of the points in ${\Cal E}$,
the error value $e_P$ may be recovered by computing 
$$e_P = {g(P)\over f'(P)},$$
where $f'$ is the formal derivative.
\endproclaim

\bigskip
Note that this is exactly the way error values are usually
determined in Reed-Solomon decoding.  The formal derivative
does not vanish at $P$ because the roots of $f$ are distinct.

\bigskip
\noindent
{\bf Proof.} We use the formula (4.4) for ${\Cal S}_e$ 
and (3.3), retaining
only terms with $e_P \ne 0$:
$${\Cal S}_e = \sum_{P \in {\Cal E}} \left(e_P \sum_{u\in \bZ_{\ge 0}^s} P^u X^{-u}\right).$$
Since $f$ is a univariate polynomial in $X_i$, nonzero terms
in the purely positive part $(f{\Cal S}_e)_+$ can only come
from terms in ${\Cal S}_e$ where the monomial $X^u$ contains
no variable other than $X_i$. (Any other terms in the product are ``mixed''
and project to zero.) As a result
$$\eqalign{
(f{\Cal S}_e)_+ &= \left(f \sum_{P\in {\Cal E}} \left(e_P \sum_{j\ge 0} X_i(P)^jX_i^{-j}\right)\right)_+\cr
 &=\left(f \sum_{P\in {\Cal E}} e_P {X_i\over X_i - X_i(P)}\right)_+ \cr
 &=X_i \sum_{P\in {\Cal E}} e_P \prod_{Q\in {\Cal E}, Q\ne P}(X_i - X_i(Q)).\cr}$$    
The polynomial $g$ appears on the right of the final line here, and
the other claims now follow from the usual analysis in the univariate
case or (3.4). $\square$

\bigskip
The same reasoning shows that in case ${\Cal E}$ is not 
in general position with respect to $X_i$ and $f$ is 
the generator for $I_{\Cal E} \cap \bF_q[X_i]$, then we still have
$(f{\Cal S})_+ = X_ig$ for $g \in \bF_q[X_i]$, but now for each
root $a$ of $f(X_i) = 0$,
$${g(a)\over f'(a)} = \sum_{P\in {\Cal E}, X_i(P) = a} e_P.$$

Even when ${\Cal E}$ is not in general position with respect
to any of the variables, the error values $e_P$ can 
be recovered from ${\Cal S}_e$ and the univariate
polynomials $f_i(X_i)$, $i = 1,\ldots,s$ generating
the collection of elimination ideals $I_{\Cal E} \cap \bF_q[X_i]$.
We illustrate the idea in a simple example with $s = 2$ before
giving the general statement.

\bigskip
\noindent
{\bf (5.2) Example.}  Let 
$${\Cal E} =\{P_1,\ldots,P_4\} = \{(0,0),(0,1),(1,1),(\alpha,1)\}$$
in ${\Bbb A}^2$ over $\bF_q$, where $\alpha \ne 0,1$.  Note that 
${\Cal E}$ is not in general position with respect to either $X$ 
or $Y$.  We have univariates $f_1(X) = X(X-1)(X-\alpha)$
and $f_2(Y) = Y(Y-1)$.  Using (3.3), (4.4), and computations as
in Example (4.6), we have
$${\Cal S}_e = e_{P_1} + e_{P_2}{Y\over Y - 1} + e_{P_3}{XY\over (X-1)(Y-1)}
+e_{P_4}{XY\over (X-\alpha)(Y-1)}.$$
(Note the special form of $h_P$ when one coordinate is zero.) Hence 
$$\eqalign{
(f_1f_2{\Cal S})_+ &= XY(e_{P_1}(X-1)(X-\alpha)(Y-1) + e_{P_2}(X-1)(X-\alpha)Y\cr
                   &\quad\ e_{P_3}X(X-\alpha)Y + e_{P_4}X(X-1)Y).\cr}\leqno(5.3)$$
Write $g(X,Y)$ for the factor in the parentheses on the right.
Note that if we substitute the points of ${\Cal E}$ in to $g$, only
one term is nonzero each time, and this allows us to determine the 
$e_{P_i}$:
$$\eqalign{g(P_1) = g(0,0) &= -\alpha e_{P_1}\cr
           g(P_2) = g(0,1) &= \alpha e_{P_2}\cr
           g(P_3) = g(1,1) &= (1-\alpha) e_{P_3}\cr
           g(P_4) = g(\alpha,1) &= \alpha(\alpha-1) e_{P_4},\cr}$$
because the factor multiplying $e_{P_i}$ is the 
product 
$$\prod_{\{\gamma : f_1(\gamma)=0,\gamma\ne X(P_i)\}} (X(P_i) - \gamma)
\prod_{\{\delta: f_2(\delta)=0,\delta\ne Y(P_i)\}} (Y(P_i) - \delta)\ne 0.\leqno(5.4)$$
There is another useful expression for (5.4).  This product is the 
same as 
$$f_1'(X(P_i))f_2'(Y(P_i)).$$
Note also that if we divide the term multiplying $e_{P_i}$ in $g(X,Y)$
by (5.4) we get one of the polynomials in a multivariable {\it Lagrange interpolation basis}
for $\bF_q[X,Y]/I_{\Cal E}$, that is a collection of polynomials
satisfying $g_i(P_k) = 0$ if $k \ne i$, and $g_i(P_k) = 1$ if 
$k = i$.  The same is true in general as we will now show.

\bigskip
\proclaim{(5.5) Proposition} Let ${\Cal E} = \{P_1, \ldots, P_t\}$
be a finite set in ${\Bbb A}^s$ over $\bF_q$.
Let $f_i$ be the monic generator of $I_{\Cal E} \cap \bF_q[X_i]$,
$i = 1,\ldots,s$.  Then 
$$(f_1f_2\cdots f_s h_{P_i})_+ = (\prod_{i=1}^s X_i) g_i,$$
where the polynomials $g_i$ satisfy
$$g_i(P_i) = \prod_{\ell = 1}^s f_\ell'(P_i),$$  
and $g_i(P_k) = 0$ if $k \ne i$.  As a result, the 
$g_i(X)/g_i(P_i)$ form a Lagrange interpolation basis for 
$\bF_q[X_1,\ldots,X_s]/I_{\Cal E}$.
\endproclaim

\bigskip
\noindent
{\bf Proof.}  This follows immediately from part (4) of (3.3)
and (3.4). $\square$

\bigskip
From (5.5) we have 
$$(f_1f_2\cdots f_s {\Cal S}_e)_+ = \left(\prod_{i=1}^t X_i\right)g$$
where $g = \sum_{i=1}^t e_{P_i} g_i.$
Hence $g(P_i) = e_{P_i} g_i(P_i)$, so by (5.5),
$$e_{P_i} = {g(P_i)\over \prod_{\ell=1}^s f_\ell'(P_i)},
\leqno(5.6)$$
and this allows us to determine the error values.

We close this section with a comment about the problem
of determining the univariate error locator polynomials $f_i$.  
This can be done easily given any Gr\"obner basis ${\Cal G}$
of $I_{\Cal E}$ (for instance the output of the BMS algorithm), using
the linear algebra techniques in $\bF_q[X_1,\ldots,X_s]/I_{\Cal E}$
described, for instance, in \cite{CLO}, Chapter 2, Section 2.
Using normal form calculations with respect to ${\Cal G}$, 
to determine $f_i$, we would simply determine the smallest
$k$ for which the normal forms of $1,X_i,X_i^2,\ldots, X_i^k$
give a linearly dependent set in $\bF_q[X_1,\ldots,X_s]/I_{\Cal E}$.
The corresponding dependence equation gives the univariate
polynomial $f_i$.  Computations of this type would 
also be used, for instance, to convert the Gr\"obner
basis ${\Cal G}$ to a lexicographic Gr\"obner basis
via the FGLM algorithm
to solve for the error locations by elimination.

\heading \S 6. Two Examples \endheading

In this section we will present two examples illustrating
the results of the previous sections.

\bigskip
\noindent
{\bf (6.1) Example.}  For our first example, we
consider Hermitian codes, in particular codes constructed
from the order domain
$R = \bF_{16}[X,Y]/\langle X^5 + Y^4 + Y\rangle$,
the affine coordinate ring of the Hermitian curve over 
$\bF_{16}$.  In the set-up from \S 2, we have $r = 1$, 
$\rho(X) = 4$, $\rho(Y) = 5$, and $\Gamma = \langle 4,5\rangle 
\subset \bZ_{\ge 0}$.  Taking the $>_{(4,5),lex}$ monomial
order the monomials in $\Delta = \{X^i Y^j : 0 \le i \le 4,
j \ge 0\}$ are an $\bF_{16}$-basis for $R$.  As is well-known,
there are 64 affine $\bF_{16}$-rational points on the Hermitian
curve $X_R$.

By the Feng-Rao bound, the minimum distance of the $C_{21} 
= ev(V_{21})^\perp$ code is at least 15, so we expect to 
be able to correct any 7 errors in a received word.  
In the order defined previously, 
$$V_{21} = Span \{1,X,Y,X^2,\ldots,Y^5\}.$$
Hence all syndromes $E_{(i,j)}$ with $i+j \le 4$ and 
$$E_{(4,1)}, E_{(3,2)}, E_{(2,3)}, E_{(1,4)}, E_{(0,5)}$$
are known initially from the received word.  In addition,
using the equation of the curve, we determine
$E_{(5,0)} = E_{(0,4)}+E_{(0,1)}$, and 
$E_{(6,0)} = E_{(1,4)}+E_{(1,1)}$.

To normalize the field, we take 
$\bF_{16} = \bF_2[\beta]/\langle \beta^4 + \beta + 1\rangle$, so
$\beta$ is a primitive element.  We consider the error of weight
7 for which 
$${\Cal E} = \{(\beta,\beta^6),(\beta^2, \beta^{14}), (\beta^4, \beta^6), (\beta^5, \beta^{14}), (\beta^8, \beta^3), (\beta^{11}, \beta^{12}), (0, 0)\}$$
and the corresponding error values are
$$\beta,\beta^4,\beta^{12},1,1,\beta,\beta.$$

As is usual for these codes, the known syndromes do not suffice
to determine the error locations and values.
Running the BMS algorithm with Feng-Rao majority voting,
additional syndromes $E_{(4,2)},E_{(3,3)},E_{(2,4)}$ 
are computed, and the curve equation furnishes the values
of $E_{(6,1)}$ and $E_{(7,0)}$.  The output of the BMS
algorithm is the following Gr\"obner basis for $I_{\Cal E}$:
$$\eqalign{ 
p_1&=X^2Y+(\beta^3+\beta^2+\beta)X+(\beta^2+\beta)Y+(\beta+1)XY+(\beta^3+\beta^2)Y^2\cr
&\qquad+(\beta^3+\beta^2+1)X^3,\cr  
p_2 &= XY^2+(\beta^2+1+\beta)Y+(\beta^3+\beta^2+\beta)Y+(\beta^2+1)X^2\cr
&\qquad+(\beta^3+\beta)XY+(\beta^2+1+\beta)Y^2+(\beta^3+\beta^2+\beta)X^3,\cr p_3&=Y^3+(\beta^3+\beta+1)X+(\beta^2+\beta)Y+\beta^2X^2+\beta^3XY\cr
&\qquad+(\beta^3+\beta^2+\beta)Y^2+(\beta^3+\beta^2+1)X^3,\cr 
p_4 &= X^4 + (\beta+1)XY+(\beta^3+\beta^2)X^3+\beta^2X^2+\beta^3 Y^2 + \beta^2Y.\cr}
\leqno(6.2)$$
The leading terms are written first in each case, so the ``footprint''
of the ideal $I_{\Cal E}$ (the set of monomials in the 
complement of $LT_>(I_{\Cal E})$) is 
$$\Delta_{\Cal E} =\{1,X,Y,X^2,XY,Y^2,X^3\},$$
and consists of the first 7 monomials in $\bF_{16}[X,Y]$ in the 
$>_{(4,5),lex}$ order.  This is the ``generic'' case for errors
of weight exactly 7 with this ordering.

At this point if we write the polynomials in (6.2) 
as $f = \sum_{m,n} X^mY^n$, then
all solve a system of equations of the form in (4.4):
$$\sum_{m,n} f_{m,n} E_{(m+r,n+s)} = 0\leqno(6.3)$$
for all 
$$(r,s)\in \{(0,0),(1,0),(0,1),(2,0),(1,1),(0,2),(3,0)\}.$$
Hence they are solutions of the truncated key equation
$$f \cdot {\Cal S}_e \equiv 0 \bmod W\leqno(6.4)$$
where $W = Span \{X^{-a} Y^{-b} : (a,b) >_{(4,5),lex} (3,0)\}$.
The polynomials in (6.2) could also be found of course by 
directly solving the linear equations (6.3).  If $W$ in (6.4)
is replaced by any $W' \subset W$, the set of solutions will
be the same.  

To determine the error values in a systematic way, we could now
exhaustively search for solutions of the system
$p_1 = \cdots = p_4 = 0$, or proceed as follows:
\roster
\item Convert the Gr\"obner basis $\{p_1,p_2,p_3,p_4\}$
to a $lex$ Gr\"obner basis using the $FGLM$ algorithm.
(If we made $X$ the smallest variable we would note that 
${\Cal E}$ is in general position with respect to $X$ as
in \S 5; with $Y$ as the smallest variable, we would see that
${\Cal E}$ is not in general position with respect to $Y$.)
\item Solve the corresponding system to find the points in 
${\Cal E}$.
\item Form any other needed univariate polynomials
in $I_{\Cal E}$ from the solutions.
\endroster

Then (5.1) with the univariate polynomial in $X$, or (5.6) will
recover the error values.  
\bigskip

One of the important things to realize about the results in 
this article is that even though this first example was constructed
using a code from an order domain with $r = 1$ (a well-studied
example of a geometric Goppa code from a curve), the actual
process of applying the BMS algorithm and determining the 
error values would be {\it exactly the same for any other example
of a $C_\ell$ code}.
This is the real lesson of \cite{HS} (although the real power
of that approach was probably not noticed at the time because
order domains of arbitrary transcendence degree had not been
used to construct codes as of yet).  At the fundamental 
level, we are always working with the ideal $I_{\Cal E}$ of a finite
set of points in ${\Bbb A}^s$, and the determination of error
locations and values can be performed in a totally uniform fashion.
 
For example, here is the same sort of computation for 
a two-dimensional extended cyclic code. (This is the dual of
the extended code corresponding to one of Hansen's toric codes, see 
\cite{H}.)  

\bigskip
\noindent
{\bf (6.4) Example.}  Let $\bF_8 = \bF_2[\alpha]/\langle \alpha^3+\alpha+1\rangle$,
and consider order domain structure on $R = \bF_8[X,Y]$ induced
by the graded lexicographic order with $X > Y$.  We have 
$$V_{10} = Span\{1,Y,X,Y^2,XY,X^2,Y^3,XY^2,X^2Y,X^3\}$$
and these give the known syndromes for $C_{10} = Ev_{10}^\perp$
(where the evaluation code is formed using all $64$ $\bF_8$
rational points in ${\Bbb A}^2$). By the Feng-Rao theorem,
this code has $d\ge 5$, so we consider an error vector
with ${\Cal E} = \{P,Q\} = \{(1,1),(\alpha,\alpha^2)\}$ and 
$e_P = 1, e_Q = \alpha^2 + 1$.

In this case the known syndromes are sufficient to determine
a Gr\"obner basis for $I_{\Cal E}$ by BMS; we are in effect
solving the truncated key equation 
$$f \cdot \overline{{\Cal S}_e} \equiv 0 \bmod W,$$
where $\overline{{\Cal S}_e}$ is the known part of the 
syndrome series, and $W = Span\{X^{-m}Y^{-n} : m+n \ge 2\}$.
The output is 
$$\{x+(\alpha^2 + \alpha) y + \alpha^2+\alpha+1, y^2 + (\alpha^2+\alpha)y+\alpha^2\}$$
which is the graded lex Gr\"obner basis for $I_{\Cal E}$.
The error values are determined using (5.1) or (5.6).
       
\refstyle{A}
 
\enddocument